\documentclass[smallextended]{svjour3}       
\smartqed  
\usepackage{graphicx}
\usepackage{cite}
\usepackage{url}

\usepackage{amsfonts}
\usepackage{amsmath} 
\usepackage{amssymb}
\usepackage{enumerate}
\usepackage{mathrsfs}
\usepackage{epstopdf}	
\usepackage{multirow}	
\usepackage{makecell}

\newcommand{\beq}{\begin{equation}}
\newcommand{\eeq}{\end{equation}}
\newcommand{\beqs}{\begin{equation*}}
\newcommand{\eeqs}{\end{equation*}}
\newcommand{\balign}{\begin{align}}
\newcommand{\ealign}{\end{align}}
\newcommand{\bea}{\begin{eqnarray}}
\newcommand{\eea}{\end{eqnarray}}
\newcommand{\ba}{\begin{array}}
\newcommand{\ea}{\end{array}}
\newcommand{\bcenter}{\begin{center}}
\newcommand{\ecenter}{\end{center}}
\newcommand{\bitem}{\begin{itemize}}
\newcommand{\eitem}{\end{itemize}}
\newcommand{\benum}{\begin{enumerate}}
\newcommand{\eenum}{\end{enumerate}}
\newcommand{\bdoc}{\begin{document}}
\newcommand{\edoc}{\end{document}}
\newcommand{\defeq}{{\stackrel{\Delta}{=}}}

\newcommand{\eg}{\mbox{\textit{e.g.}, \/}}
\newcommand{\etal}{\mbox{\textit{et al.} \/}}
\newcommand{\viz}{\mbox{\textit{\iet viz.},\ \/}}
\newcommand{\ie}{\mbox{\textit{i.e.},\ \/}}
\newcommand{\etc}{\mbox{\textit{etc.},\ \/}}
\newcommand{\cf}{\mbox{\textit{cf.},\ \/}}

\usepackage[usenames,dvipsnames,svgnames,table]{xcolor}

\newcommand{\black}[1]{\color{black} #1}
\newcommand{\yellow}[1]{\color{yellow} #1}
\newcommand{\green}[1]{\color{green} #1}
\newcommand{\bgrd}[1]{\color{bgrd} #1}
\newcommand{\magenta}[1]{\color{magenta} #1}
\newcommand{\white}[1]{\color{white} #1}
\newcommand{\gray}[1]{\color{gray} #1}

\newcommand\Quote[1]{\lq\textsl{#1}\rq}
\newcommand\fr[2]{{\textstyle\frac{#1}{#2}}}

\newcommand{\Ambb}{\mathbb{A}}
\newcommand{\Bmbb}{\mathbb{B}}
\newcommand{\Cmbb}{\mathbb{C}}
\newcommand{\Dmbb}{\mathbb{D}}
\newcommand{\Embb}{\mathbb{E}}
\newcommand{\Fmbb}{\mathbb{F}}
\newcommand{\Gmbb}{\mathbb{G}}
\newcommand{\Hmbb}{\mathbb{H}}
\newcommand{\Imbb}{\mathbb{I}}
\newcommand{\Jmbb}{\mathbb{J}}
\newcommand{\Kmbb}{\mathbb{K}}
\newcommand{\Lmbb}{\mathbb{L}}
\newcommand{\Mmbb}{\mathbb{M}}
\newcommand{\Nmbb}{\mathbb{N}}
\newcommand{\Ombb}{\mathbb{O}}
\newcommand{\Pmbb}{\mathbb{P}}
\newcommand{\Qmbb}{\mathbb{Q}}
\newcommand{\Rmbb}{\mathbb{R}}
\newcommand{\Smbb}{\mathbb{S}}
\newcommand{\Tmbb}{\mathbb{T}}
\newcommand{\Umbb}{\mathbb{U}}
\newcommand{\Vmbb}{\mathbb{V}}
\newcommand{\Wmbb}{\mathbb{W}}
\newcommand{\Xmbb}{\mathbb{X}}
\newcommand{\Ymbb}{\mathbb{Y}}
\newcommand{\Zmbb}{\mathbb{Z}}

\newcommand{\Amsc}{\mathcal{A}}
\newcommand{\Bmsc}{\mathcal{B}}
\newcommand{\Cmsc}{\mathcal{C}}
\newcommand{\Dmsc}{\mathcal{D}}
\newcommand{\Emsc}{\mathcal{E}}
\newcommand{\Fmsc}{\mathcal{F}}
\newcommand{\Gmsc}{\mathcal{G}}
\newcommand{\Hmsc}{\mathcal{H}}
\newcommand{\Imsc}{\mathcal{I}}
\newcommand{\Jmsc}{\mathcal{J}}
\newcommand{\Kmsc}{\mathcal{K}}
\newcommand{\Lmsc}{\mathcal{L}}
\newcommand{\Mmsc}{\mathcal{M}}
\newcommand{\Nmsc}{\mathcal{N}}
\newcommand{\Omsc}{\mathcal{O}}
\newcommand{\Pmsc}{\mathcal{P}}
\newcommand{\Qmsc}{\mathcal{Q}}
\newcommand{\Rmsc}{\mathcal{R}}
\newcommand{\Smsc}{\mathcal{S}}
\newcommand{\Tmsc}{\mathcal{T}}
\newcommand{\Umsc}{\mathcal{U}}
\newcommand{\Vmsc}{\mathcal{V}}
\newcommand{\Wmsc}{\mathcal{W}}
\newcommand{\Xmsc}{\mathcal{X}}
\newcommand{\Ymsc}{\mathcal{Y}}
\newcommand{\Zmsc}{\mathcal{Z}}

\newcommand{\Ascr}{\mathscr{A}}
\newcommand{\Bscr}{\mathscr{B}}
\newcommand{\Cscr}{\mathscr{C}}
\newcommand{\Dscr}{\mathscr{D}}
\newcommand{\Escr}{\mathscr{E}}
\newcommand{\Fscr}{\mathscr{F}}
\newcommand{\Gscr}{\mathscr{G}}
\newcommand{\Hscr}{\mathscr{H}}
\newcommand{\Iscr}{\mathscr{I}}
\newcommand{\Jscr}{\mathscr{J}}
\newcommand{\Kscr}{\mathscr{K}}
\newcommand{\Lscr}{\mathscr{L}}
\newcommand{\Mscr}{\mathscr{M}}
\newcommand{\Nscr}{\mathscr{N}}
\newcommand{\Oscr}{\mathscr{O}}
\newcommand{\Pscr}{\mathscr{P}}
\newcommand{\Qscr}{\mathscr{Q}}
\newcommand{\Rscr}{\mathscr{R}}
\newcommand{\Sscr}{\mathscr{S}}
\newcommand{\Tscr}{\mathscr{T}}
\newcommand{\Uscr}{\mathscr{U}}
\newcommand{\Vscr}{\mathscr{V}}
\newcommand{\Wscr}{\mathscr{W}}
\newcommand{\Xscr}{\mathscr{X}}
\newcommand{\Yscr}{\mathscr{Y}}
\newcommand{\Zscr}{\mathscr{Z}}

\newcommand\independent{\protect\mathpalette{\protect\independenT}{\perp}}
\def\independenT#1#2{\mathrel{\rlap{$#1#2$}\mkern2mu{#1#2}}}

\newcommand{\bsf}[1]{\textsf{\textbf{#1}}}
\newcommand{\lbsf}[1]{\textsf{\large  \textbf{#1}}}
\newcommand{\Lbsf}[1]{\textsf{\Large  \textbf{#1}}}
\newcommand{\hbsf}[1]{\textsf{\huge  \textbf{#1}}}

\newcommand{\myminipage}[3]{\begin{minipage}[#1]{#2}{#3} \end{minipage}}
\newcommand{\sbs}[4]{\myminipage{c}{#1}{#3} \hfill
\myminipage{c}{#2}{#4}}

\newcommand{\myfig}[2]{\centerline{\psfig{figure=#1,width=#2,silent=}}}
\newcommand{\myfigh}[2]{\centerline{\psfig{figure=#1,height=#2,silent=}}}
\newcommand{\myfigwh}[3]{\centerline{\psfig{figure=#1,width=#2,height=#3,silent=}}}

\newcommand{\beqa}{\begin{eqnarray}}
\newcommand{\eeqa}{\end{eqnarray}}
\newcommand{\beqan}{\begin{eqnarray*}}
\newcommand{\eeqan}{\end{eqnarray*}}
\newcommand{\dst}[1]{\displaystyle{ #1 }}
\newcounter{pf}

\newcommand{\smax}[1] { \bar \sigma \left( #1 \right) }
\newcommand{\Rn}{{\mathbb R}^n}
\newcommand{\R}{{\mathbb R}}
\newcommand{\C}{{\mathbb C}}
\newcommand{\Rm}{\mathbb{R}^m}
\newcommand{\Rmn}{\mathbb{R}^{m \times n}}
\newcommand{\Rpq}{\mathbb{R}^{p \times q}}
\newcommand{\Cn}{\mathbb{C}^n}
\newcommand{\Cm}{\mathbb{C}^m}
\newcommand{\Cnn}{\mathbb{C}^{n \times n}}
\newcommand{\Cmn}{\mathbb{C}^{m \times n}}
\newcommand{\ip}[1]{\left\langle #1 \right\rangle}
\newcommand{\rank}{\mbox{rank}}
\newcommand{\Span}{\mbox{\rm Span }}
\newcommand{\Trace}{\mbox{\rm Tr }}
\newcommand{\Spec}{\mbox{\rm Spec }}
\newcommand{\vectornorm}[1]{\left\|#1\right\|}

\newcommand{\pd}[2]{\frac{\partial #1}{\partial #2}}
\newcommand{\ppd}[3]{\frac{\partial^2 #1}{\partial #2 \partial #3}}

\newcommand{\thtilde}{\tilde{\theta}}
\newcommand{\thnom}{\theta^\circ}
\newcommand{\thopt}{\theta^{\mbox{\small opt}}}
\newcommand{\thhat}{{\hat{\theta}}}
\newcommand{\Tho}{\Theta^\circ}
\newcommand{\tho}{\theta^\circ}
\newcommand{\np}{{n_p}}

\newcommand{\ii}{{[i]}}
\newcommand{\II}{{[i+1]}}
\newcommand{\iii}{{[ii]}}
\newcommand{\jj}{{[j]}}
\newcommand{\kk}{{[k]}}
\newcommand{\thi}{{\theta^\ii}}
\newcommand{\thI}{{\theta^\II}}
\newcommand{\di}{{d^\ii}}
\newcommand{\gi}{{g^\ii}}
\newcommand{\Hi}{{\HH^\ii}}
\newcommand{\thK}{\theta^{(k+1)}}
\newcommand{\gk}{{g^{(k)}}}
\newcommand{\Hk}{{{\cal H}^{(k)}}}

\newcommand{\bfdelta}{{\bf \Delta}}

\newcommand{\Exp}[1]{\exp \left\{ #1 \right\}} 
\newcommand{\gaussian}[1]{\mathbb{N} \left( #1 \right)}
\newcommand{\uniform}[1]{\mathbb{U} \left[ #1 \right]}
\newcommand{\exponential}[1]{\mathbb{E} \left[ #1 \right]}
\newcommand{\EXP}[1]{\EEXP \left[ #1 \right]} 
\newcommand{\EEXP}{\mbox{\bsf{E}}} 
\newcommand{\Prob}[1]{\mbox{{\sf Pr}} \left(#1 \right)}
\newcommand{\convas}{\stackrel{as}{\longrightarrow}}
\newcommand{\convinp}{\stackrel{p}{\longrightarrow}}
\newcommand{\convind}{\stackrel{d}{\longrightarrow}}
\newcommand{\convqm}{\stackrel{qm}{\longrightarrow}}
\newcommand{\sss}[1]{{_{#1}}}
\newcommand{\density}[2]{p_{_{_{#1}}}\!\!\left(#2 \right)} 
\newcommand{\distro}[2]{P_{_{_{#1}}}\!\!\left(#2 \right)} 
\newcommand{\rxx}[1]{R_{_{#1}}\!} 
\newcommand{\sxx}[1]{S_{_{#1}}} 
\newcommand{\cov}[1]{\Lambda_{_{#1}}} 
\newcommand{\mean}[1]{m_{_{#1}}} 
\newcommand{\LS}[1]{\hat{#1}_{_{LS}}} 
\newcommand{\MV}[1]{\hat{#1}_{_{MV}}} 
\newcommand{\LMV}[1]{\hat{#1}_{_{LMV}}} 
\newcommand{\ML}[1]{\hat{#1}_{_{ML}}} 

\renewcommand{\arraystretch}{0.9}
\newcommand{\bmat}[1]{ \begin{bmatrix} #1 \end{bmatrix}}
\newcommand{\mat}[1]{ \left[ \begin{array}{cccccccc} #1 \end{array}
\right] }
\newcommand{\smallmat}[1]{\small{\mat{#1}}}
\newcommand{\sysblk}[4]{\begin{array}{c|cccc}#1&#2\\ \hline#3&#4
\end{array}}
\newcommand{\sysmat}[4]{\left[\sysblk{#1}{#2}{#3}{#4}\right]}
\newcommand{\SGeq}{\succ}
\newcommand{\SLeq}{\prec}
\newcommand{\Geq}{\succeq}
\newcommand{\Leq}{\preceq}

\newcommand{\Bset}{\mathbb{B}}
\newcommand{\Cset}{\mathbb{C}}
\newcommand{\Fset}{\mathbb{F}}
\newcommand{\Mset}{\mathbb{M}}
\newcommand{\Nset}{\mathbb{N}}
\newcommand{\Qset}{\mathbb{Q}}
\newcommand{\Rset}{\mathbb{R}}
\newcommand{\Sset}{\mathbb{S}}
\newcommand{\Tset}{\mathbb{T}}
\newcommand{\Uset}{\mathbb{U}}
\newcommand{\Vset}{\mathbb{V}}
\newcommand{\Wset}{\mathbb{W}}
\newcommand{\Zset}{\mathbb{Z}}

\newcommand{\Ical}{{\cal I}}
\newcommand{\Acal}{{\cal A}}
\newcommand{\Bcal}{{\cal B}}
\newcommand{\Ccal}{{\cal C}}
\newcommand{\Dcal}{{\cal D}}
\newcommand{\Ecal}{{\cal E}}
\newcommand{\Fcal}{{\cal F}}
\newcommand{\Gcal}{{\cal G}}
\newcommand{\Hcal}{{\cal H}}
\newcommand{\Kcal}{{\cal K}}
\newcommand{\Lcal}{{\cal L}}
\newcommand{\Mcal}{{\cal M}}
\newcommand{\Ncal}{{\cal N}}
\newcommand{\Pcal}{{\cal P}}
\newcommand{\Qcal}{{\cal Q}}
\newcommand{\Rcal}{{\cal R}}
\newcommand{\Scal}{{\cal S}}
\newcommand{\Tcal}{{\cal T}}
\newcommand{\Wcal}{{\cal W}}
\newcommand{\Ucal}{{\cal U}}
\newcommand{\Vcal}{{\cal V}}
\newcommand{\Xcal}{{\cal X}}
\newcommand{\Zcal}{{\cal Z}}

\newcommand{\FF}{{\bf F}}
\newcommand{\GG}{{\bf G}}
\newcommand{\HH}{{\bf H}}
\newcommand{\LL}{{\bf L}}
\newcommand{\NN}{{\bf N}}
\newcommand{\MM}{{\bf M}}
\newcommand{\PP}{{\bf P}}
\newcommand{\QQ}{{\bf Q}}
\newcommand{\RR}{{\bf R}}
\renewcommand{\SS}{{\bf S}}
\newcommand{\TT}{{\bf T}}
\newcommand{\VV}{{\bf V}}
\newcommand{\WW}{{\bf W}}

\newcommand{\thk}{\theta^{(k)}}
\newcommand{\thb}{\theta^{\rm opt}}
\newcommand{\alb}{\alpha^{\rm opt}}
\newcommand{\dk}{d^{(k)}}
\newcommand{\Hinf}{{\cal H}_\infty}
\newcommand{\Htwo}{{\cal H}_2}

\renewcommand{\arraystretch}{1.1}

\newcommand{\red}[1]{{\color{red} #1}}
\newcommand{\blue}[1]{{\color{Blue} #1}}

\newcounter{l1}
\newcounter{l2}
\newcounter{l3}
\newcommand{\bdotlist}{\begin{list}{$\bullet$}{}}
\newcommand{\bboxlist}{\begin{list}{$\Box$}{}}
\newcommand{\bbboxlist}{\begin{list}{\raisebox{.005in}{{\tiny
$\blacksquare$ \ \ }}}{}}
\newcommand{\bdashlist}{\begin{list}{$-$}{} }
\newcommand{\blist}{\begin{list}{}{} }
\newcommand{\barablist}{\begin{list}{\arabic{l1}}{\usecounter{l1}}}
\newcommand{\balphlist}{\begin{list}{(\alph{l2})}{\usecounter{l2}}}
\newcommand{\bAlphlist}{\begin{list}{\Alph{l2}.}{\usecounter{l2}}}
\newcommand{\bdiamlist}{\begin{list}{$\diamond$}{}}
\newcommand{\bromalist}{\begin{list}{(\roman{l3})}{\usecounter{l3}}}

\newcommand{\thm}[1]{\noindent \begin{theorem} #1   \end{theorem}}
\newcommand{\prop}[1]{\begin{proposition} #1 \end{proposition}}
\newcommand{\lem}[1]{\begin{lemma} #1  \hfill $\blacksquare$ \end{lemma}}
\newcommand{\ex}[1]{\begin{example} {\rm #1} \end{example}}
\newcommand{\prf}[1]{ \noindent {\em Proof:} \, #1 \hfill $\blacksquare$}
\newcommand{\rem}[1]{\begin{remark} {\rm #1} \hfill $\Box$ \end{remark}}
\newcommand{\defn}[1]{\begin{definition} {\rm #1 } \end{definition}}
\newcommand{\prob}[1]{\begin{exercise} {\rm  #1 } \end{exercise}}
\newcommand{\cor}[1]{\begin{corollary}   #1  \end{corollary}}

\newcommand{\argmin}{\mathop{\rm argmin}}
\newcommand{\argmax}{\mathop{\rm argmax}}
\newcommand{\diag}{\mathop{\mathrm{diag}}}
\newcommand{\tr}{\mathop{\rm Tr}}
\newcommand{\conv}{\mathop{\rm conv}}
\newcommand{\var}{\mathop{\rm Var}}
\renewcommand{\b}[1]{\ensuremath{\boldsymbol{\mathrm{#1}}}}
\newcommand{\ms}{{\rm MS}}
\newcommand{\tcs}{{\rm TCS}}
\newcommand{\scs}{{\rm SCS}}

\newcommand{\E}[1]{\b{\mu}_{{#1}}}
\newcommand{\Var}[1]{{\Sigma_{#1}}}

\def\u{u}	
\def\U{U}	
\def\v{v}	
\def\V{V}	
\def\s{s}	
\def\S{S}	
\def\f{q}	
\def\q{q}	
\def\Q{Q}	
\def\c{c}	
\def\z{z}	
\def\xcap{b}
\def\ccap{\mathbf{\c}}
\def\rentftr{\Phi}
\def\rentfsr{\Sigma}
\def\rentDer{\Delta}
\def\lambdao{\lambda^{o}}
\def\RT{{\rm \Phi}(\Tcal, \Fcal)}
\def\RS{{\rm \Sigma}(\Scal, \Ecal)}

\newtheorem{assumption}{Assumption}

\journalname{}
\begin{document}

\title{Financial Storage Rights in Electric Power Networks
}


\author{Daniel Mu\~noz-\'Alvarez$^\dagger$         \and
        Eilyan Bitar$^{\dagger\dagger}$ 
}

\authorrunning{D. Mu\~noz-\'Alvarez, E. Bitar} 

\institute{$\dagger$ School of Electrical and Computer Engineering, Cornell University,
              Ithaca, NY, 14853, \email{dm634@cornell.edu}  \\
              $\dagger\dagger$ 
              School of Electrical and Computer Engineering, and School of Operations Research and Information Engineering, Cornell University, 
              Ithaca, NY, 14853, \email{eyb5@cornell.edu} }

\date{}

\maketitle

\begin{abstract}
The decreasing cost of energy storage technologies coupled with their potential to bring significant benefits to electric power networks have kindled research efforts to design both market and regulatory frameworks to facilitate the efficient construction and operation of  such technologies.
In this paper, we examine an \emph{open access} approach to the integration of storage, which enables the complete decoupling of a storage facility's ownership structure from its  operation.
In particular, we analyze a nodal spot pricing system built on a model of economic dispatch in which storage is centrally dispatched by the independent system operator (ISO) to maximize social welfare.
Concomitant with such an approach is the ISO's collection of a merchandising surplus reflecting congestion in storage.
We introduce a class of tradable electricity derivatives  -- referred to as \emph{financial storage rights} (FSRs) -- to enable the redistribution of such rents in the form of financial property rights to storage capacity; and establish a generalized \emph{simultaneous  feasibility test} to ensure the ISO's revenue adequacy  when allocating such financial property rights to market participants. 
Several advantages of such an approach to open access storage are discussed. In particular, we illustrate with a stylized example  the role of FSRs in synthesizing fully hedged, fixed-price bilateral contracts for energy, when the seller and buyer exhibit differing intertemporal supply and demand characteristics, respectively.

\end{abstract}

\section{Introduction} \label{sec_introduction}

The increased penetration of supply derived from variable renewable energy resources, coupled with the recent decline in the cost of electric energy storage technologies, has brought about an opportunity to significantly reduce the cost of managing the electric power system through careful planning, deployment, and operation of storage resources \cite{HindsBoyer14}. Broadly, the short-run value of energy storage derives from its ability to arbitrage energy  forward in time, enabling both the absorption of power imbalances on short time scales and the more substantial reshaping of supply and demand profiles over longer periods of time.
The extent to which the deployment of a collection of energy storage devices might benefit the power system depends critically, however, on the collective \emph{sizing}, \emph{placement}, and \emph{operation} of said devices \cite{BitarBose14}. The challenge resides in the design and implementation of electricity markets and instruments that induce strategic expansion and operation of storage in a manner that  is consistent with the maximization of social welfare over both the long and short run, respectively. 

The coordinated optimal  dispatch of a collection of  distributed energy storage resources clearly offers the possibility of a sizable reduction in the cost  of servicing demand by reshaping  it in such a manner as to alleviate both transmission congestion and the reliance on peak power generation \cite{PJM12}. Of interest then is the characterization of mechanisms for the integration of storage, which encourage its efficient operation.  And of critical importance to this effort is the resolution of the question: \emph{who commands the storage?}
Among the variety of possible answers to this question, there are two extremes -- differing in terms of the degree of government intervention -- which we naturally refer to as \emph{competitive} and \emph{regulated}. Each implies a distinct mechanism for both the operation of the physical storage facilities and the remuneration of the services provided.

Broadly, the \emph{competitive} or \emph{market-based operation of storage} entails a decentralized operating paradigm in which storage owners pursue their own rational (profit maximizing) interests in the spot energy market.
A shortcoming of such approach to storage integration derives from the uncertainty in revenue that storage owner-operators might obtain from the spot market. Energy storage is a  capital intensive technology.
And several recent studies  \cite{DruryEtAl11, Sioshansietal12, DOE13, RMI15b} have indicated that the risk of incomplete capital cost recovery due to such revenue uncertainty  may serve to inhibit investment in storage facilities.
Sioshansi \cite{Sioshansi10} also goes on to show that a complete reliance on the spot energy market to guide the integration of storage may lead to its substantial underutilization relative to the social optimum, as strategic owner-operators of storage will naturally endeavor to preserve intertemporal price differences for purposes of arbitrage.

The \emph{regulated operation of storage}, on the other hand, calls for a centralized operating paradigm in which storage is treated as a communal asset that is centrally dispatched by the  Independent System Operator (ISO) to maximize social welfare subject to its physical constraints.\footnote{The PJM Interconnection has explored a similar regulatory framework in which energy storage would be operated and compensated traditionally like a transmission asset \cite{PJM12}.}
The socially optimal dispatch of storage, in concert with conventional generation and transmission, naturally improves upon the welfare of the system in the short run.
Accordingly, such an approach to the operation of storage necessitates the creation of a mechanism capable of extracting and redistributing the value added by storage back to the owners of the responsible storage facilities.
Towards this end, we propose a market mechanism founded on the definition of tradable financial instruments, which monetize property rights to storage capacity made available to the ISO for centralized operation.
Such an approach resembles the regulation and operation of transmission in the majority of US electricity markets, which entails the centrally optimized operation of the transmission network subject to the locational marginal pricing of energy, and the allocation of financial transmission rights that monetize  property rights to said transmission capacity  \cite{AlsacEtAl04, Hogan92, Hogan02, ONeillEtAl02, ONeillEtAl13}.

\subsection{Open Access Energy Storage}
 There has been  recent activity in both academia and industry  to identify alternative paradigms to support the efficient integration of storage into power system operations \cite{Sioshansietal12, PJM12}.  One stream of literature 
centers on an \textit{open access} approach to the  integration of storage; or more simply, \textit{open access storage} (OAS) \cite{Heetal11, Sioshansietal12, Taylor14b}.
Loosely, we refer to OAS as a regulatory framework in which energy storage facilities are treated as communal assets accessible by all participants in the wholesale energy market. 

To the best of our knowledge, only two concrete approaches to OAS have been proposed. He et al. \cite{Heetal11} proposes a market framework where storage owners sell \emph{physically binding rights} to their storage capacity through sequential auctions coordinated by the ISO.
The collection of physical rights, which are defined as a sequence of nodal power injections within a specified time horizon, determine the actual operation of the storage.
As such, the physical rights associated with a particular storage facility must be collectively feasible with respect to the corresponding physical device constraints. 
While such physical rights might be used by market participants to execute price arbitrage or mitigate the cost of honoring existing contractual energy commitments, there are several important limitations.
First, the ability of a market participant to leverage on a physical storage right depends on her location within the network relative to the storage facilities. Such restriction could serve to limit market access.
Second, the eventual physical dispatch of storage is determined by  a sequence of auctions -- the outcome of which is likely to substantially deviate from the socially optimal dispatch, because of strategic interactions between parties bidding for physical storage rights.

Closer to our proposal, Taylor \cite{Taylor14b} suggests an approach to OAS that centers on a paradigm in which storage owners sell \emph{financially binding rights} to their storage capacity through an auction  coordinated by the ISO.
The ISO is charged with the task of  operating storage in a socially efficient manner -- not unlike its non-discriminatory operation of the transmission network.
As financial rights, they do not interfere with the optimal operation of storage, but rather, they represent entitlements to portions of the merchandising surplus collected by the ISO.
A central component of the proposal in \cite{Taylor14b} is the definition of the financial rights in terms of the shadow prices associated with the physical constraints on the storage facilities. This  is analogous to the definition of flowgate rights (FGRs) \cite{ChaoPeck96, ChaoPeck97, Chaoetal00} in the context of open access transmission.
And, as a result, such a definition of financial storage rights is naturally endowed with advantages and disadvantages comparable to those of FGRs in the context of transmission. We refer the reader to Section \ref{sec:fsr} and \cite{Hogan00, Orenetal95} for a more detailed discussion on such issues.

\subsection{Contribution}

We propose a regulatory framework to enable open access storage, which centers largely on the concept of \textit{financial storage rights} (FSRs). Broadly speaking, FSRs can be interpreted as  financial property rights to storage capacity; or, more accurately, as financial entitlements (or obligations) to portions of the storage congestion rent collected by the ISO under the socially optimal dispatch of storage capacity. Being defined as such, FSRs enable the complete decoupling of  a storage facility's  ownership from  its physical operation.
Moreover, the specific form of FSRs that we propose -- viz. a sequence of nodal power injections and withdrawals that yield its holder a payment according to the corresponding sequence of nodal spot prices  -- provides market participants the ability to perfectly hedge physical or financial energy positions against intertemporal price risk in the spot market.\footnote{Such a definition of FSRs represents a financial analog to the physical storage rights proposed by He et al. \cite{Heetal11}, and is in contrast to the  constraint-based financial rights proposed in \cite{Taylor14b}.} Such hedging capabilities represent a natural complement to financial transmission rights (FTRs) and their ability to hedge spatial price risk across the network. 
What distinguishes such financial instruments from standard forward energy contracts is the fact that they are issued under the physical cover of transmission and storage capacity, and are settled against the merchandising surplus collected by the ISO. Accordingly, in Section  \ref{sec:generalSimultaneousFeasibility}, we establish a generalized \emph{simultaneous  feasibility test} (SFT), which  constrains the joint allocation of financial transmission and storage rights in such a manner as to guarantee the ISO's revenue adequacy. Namely, any simultaneously feasible collection of transmission and storage rights are guaranteed to yield a rent that does not exceed the merchandising surplus collected by the ISO. A positive attribute of the proposed SFT is that it enables the allocation (auction) of FSRs at nodes without physical storage capacity -- a feature which genuinely democratizes access to storage by all market participants.

\subsection{Organization}

\noindent The remainder of paper is organized as follows. 
In Section \ref{sec_Formulation}, we formulate the  multi-period economic dispatch  problem with storage, and delineate its optimality conditions. 
In Section \ref{sec_FSR}, we formally introduce the concept of financial storage rights, and establish a general test for simultaneous  feasibility, which restricts the allocation of both financial transmission and storage rights in such a manner as to ensure the ISO's revenue adequacy.
In Section \ref{sec_example}, we illustrate with a stylized example the role of FSRs in synthesizing flexible, fully hedged, fixed-price bilateral contracts for energy.
We close with a discussion on directions for future research in Section \ref{sec_Conclusions}.
All mathematical proofs  are included in the Appendix to the paper.

\section{Models and Formulation} \label{sec_Formulation}

\subsection{Notation}

\noindent Let $\Rset$ denote the set of real numbers and $\Rset_+$ the non-negative real numbers. Denote the transpose of a vector $\b{x} \in \Rset^n$ by $\b{x}^{\top}$.
Let $x_i$ denote the $i^{\rm th}$ entry of a vector $\b{x} \in \Rset^n$.
We define by $\b{1}$  a column vector of all ones and by $\b{e}_i$ the $i^{\rm th}$ standard basis vector of dimension appropriate to the context. For two matrices $A,B \in \Rset^{m\times n}$ of equivalent dimension, we denote their Hadamard product by $A \circ B$. Given a matrix $A \in  \Rset^{m\times n}$, we write $A = 0$ to denote entrywise equivalence to zero.

\subsection{Network Model} \label{sec:network_mod}

\noindent Consider a transmission network defined on a set  of $n$ nodes (buses) connected by $m$ edges (transmission lines). The associated graph of the network is assumed connected. The nodes are indexed by $i = 1,2, \dots,n$. We operate under the assumption of a linear  model of steady state power flow defined by the so called DC power flow approximation, where the vector of nodal power injections is linearly mapped to a vector  of (directional) power flows along the $m$ transmission lines through the mapping  $H \in \Rset^{2m \times n}$, commonly referred to as the \emph{shift-factor matrix}. Let $\b{c} \in \Rset^{2m}_+$ denote the corresponding vector of transmission line capacities. It follows that the set of feasible power injections is described by the polytope $\Pcal(\b{c}) \subset \Rset^n$,
\begin{align} \label{eq_Polytope}
\Pcal(\b{c}) \ = \ \left\{ \b{v} \in \Rmbb^{n} \ \left| \ H\b{v} \leq \b{c}, \  \ \b{1}^{\top}\b{v} = 0    \right.\right\}.
\end{align}
One can readily verify the compactness of $\Pcal(\b{c})$, as $\rank(H) = n-1$ and $\b{1}^{\top}$ is linearly independent from the rows of $H$.

\subsection{Cost Model}
\noindent At the core of the formulation considered in this paper is the problem of multi-period economic dispatch over $N$ discrete time periods, which we index by $k=0,\ldots,N-1$. 
We measure the cost and benefit of the net injection vector $\b{v}(k) \in \Rset^n$ at time $k$ according to
$$ C(\b{v}(k), k) = \sum_{i=1}^n C_i(v_i(k),k), $$
Each component function $C_i(v,k)$ is assumed to be increasing, convex, and  differentiable in $v$ over $\Rset$. Moreover, each function is assumed to satisfy $C_i(0,k) = 0$, \ $C_i(v,k) > 0$ for $v > 0$, and $C_i(v,k) <  0$ for $v< 0$. 
This implies that $C_i(v,k)$ represents the \emph{convex cost of generation} for $v>0$ at node $i$ and time $k$. Conversely,  $-C_i(v,k)$ represents the \emph{concave benefit of consumption} for $v<0$ at node $i$ and time $k$. Finally, the component functions $\{C_i(\cdot,k)\}$ to allowed to vary with time in order to capture the potential variation in the nodal demand preferences over time. We refer the reader to Wu et al. \cite{Wuetal96} for a more detailed explanation of this model. 

\subsection{Energy Storage Model} \label{sec:storage model}
\noindent We consider an arbitrary collection of $n$ perfectly efficient energy storage devices connected to the transmission network, where we associate with each node $i$ a storage device with energy capacity $b_i  \in \Rset_+$.  We denote by $\b{b} = \bmat{b_1,\dots, b_n}^{\top}$ the vector of nodal energy storage capacities. 
The collective storage dynamics are naturally modeled as a linear difference equation
\beq
\label{eq:storagemod} \b{z}(k+1) = \b{z}(k)  - \b{u}(k)  \ 
\eeq
for $k=0,1,\dots,N-1$, where the vector $\b{z}(k) \in \Rset_+^n$ denotes the vector of energy storage states just preceding time period $k$, and the input $\b{u}(k) \in \Rset^n$  denotes the vector of net energy storage extractions during period $k$. The notational convention is such that $u_i(k) > 0$ (resp. $u_i(k) < 0$) represents a net energy extraction from (resp. injection into) the  storage device at node $i$ during time period $k$. 
Without loss of generality, we assume an initial condition  of $\b{z}(0) = 0$ for the remainder of the paper. 
The limited capacities of the energy storage devices require that $ 0 \leq \b{z}(k) \leq \b{b}$ for all $k$.
Iterating the linear difference equation (\ref{eq:storagemod}) back to its initial condition, one can express the storage capacity constraint as
\begin{align}
0 \leq  \  -\sum_{\ell =0}^{k-1} \b{u}(\ell)  \  \leq \b{b}  \label{eq:con1}
\end{align}
for $k = 1,\dots, N$.
As a matter of notational convenience, we consider an \textit{equivalent} characterization of the energy storage capacity constraints (\ref{eq:con1}), which enables a decomposition of the constraints across nodes. More specifically, 
letting $\b{u}_i = \bmat{u_i(0), \dots, u_i(N-1)}^{\top}$ denote the entire sequence of injections and extractions from the storage device at node $i$,  one can recast the constraints defined by (\ref{eq:con1}) as
\begin{align}
\b{u}_i \ \in \ \Ucal(\b{b}_i)   \ = \ \left\{  \left. \b{u} \in \Rset^N \ \right| \  0 \leq L \b{u} \leq  \b{b}_i   \right\}               \label{stg_model}
\end{align}
for $i = 1,\dots, n$. Here,  $L \in \Rmbb^{N\times N}$ denotes a lower triangular matrix with entries $[L]_{k \ell} = -1$ for all $ k \ge \ell$, and zero otherwise. We also define  $\b{b}_i = \bmat{b_i, \dots, b_i}^\top \in \Rset^n$.
It is immediate to see that $\Ucal(\b{b}_i)$ is a compact polytope containing the origin for each $i = 1,\dots,n$.

\begin{remark} 
While the model of storage considered is stylized in nature, much of the ensuing analysis and conclusions derived can be easily extended to accommodate nonidealities  in storage, such as constraints on allowable rates of charging and discharging, roundtrip inefficiencies, and dissipative losses.
\end{remark}

\subsection{Multi-Period Economic Dispatch}
Working within the idealized setting considered, 
we now formulate the problem of multi-period economic dispatch with storage.
Broadly, the objective of the ISO is to select a vector of nodal prices for energy that sustains a competitive equilibrium between supply and demand at a feasible system operating point that maximizes social welfare -- a so-called \emph{economic dispatch}. 
Formally, the \emph{multi-period economic dispatch} problem is stated as:
\begin{alignat}{3}   
& \text{minimize}  \hspace{.25in} & &	\sum_{k=0}^{N-1} \ C(\b{v}(k),k) \label{objective} \\
& \text{subject to} 	& &	 \b{v}(k)  + \b{u}(k)  \in  \  \Pmsc(\b{c}),  \qquad  k=0,\dots,N-1   \label{constraints_flow}  \\
					&& & 	\hspace{.58in} \b{u}_{i}  \in  \  \Umsc( \b{b}_i),  \hspace{.2in} i =1,\dots,n \label{constraints_stg}  
\end{alignat} 
where the minimization is taken with respect to the variables $\b{v}(k) \in \Rset^n$ and  $\b{u}(k) \in \Rset^n$ for $k =0,\dots,N-1$. We will occasionally denote the decision variables more compactly by the pair $(V,U)$, where $V = \bmat{\b{v}(0), \dots, \b{v}(N-1)}$ and $U = \bmat{\b{u}(0), \dots, \b{u}(N-1)}$.

\subsection{Optimality Conditions}

\begin{definition}
\label{def:feas}
A pair $(V,U)$ is a \emph{feasible dispatch}  if it satisfies constraints \eqref{constraints_flow}-\eqref{constraints_stg}. A pair $(V,U)$ is an (optimal) \emph{economic dispatch} if it solves problem (\ref{objective})-(\ref{constraints_stg}).
\end{definition}

The multi-period economic dispatch problem (\ref{objective}) - (\ref{constraints_stg}) is a convex optimization problem with linear constraints.
As such, an economic dispatch $(V, U)$ is characterized by the existence of Lagrange multipliers such that the Karush-Kuhn-Tucker (KKT) conditions (\ref{constraints_flow}) - (\ref{FOC_7}) hold.
More specifically, we associate Lagrange  multipliers $\gamma(k) \in \Rset$  and $\b{\mu}(k) \in \Rset_+^{2m}$ with the power balance and line flow capacity constraints (\ref{constraints_flow}) at time $k$, respectively. Similarly, we define  $\underline{\b{\nu}}_i \in \Rset_+^N$ and $\overline{\b{\nu}}_i \in \Rset_+^N$  as the Lagrange multipliers associated with the energy capacity constraints (\ref{constraints_stg}) of the storage device at node $i$. 
In specifying  the KKT conditions, it will be convenient to define as $\b{\lambda}(k) \in  \Rset^n$ a particular linear combination of Lagrange multipliers given by:
\begin{align}  \label{eq:LMP}
\b{\lambda}(k) = \gamma(k) \b{1}  - H^{\top} \b{\mu}(k)
\end{align}
for each time $k=0, \dots, N-1$.  The \emph{stationarity condition} is given by:
\begin{align}
\nabla  C(\b{v}(k), k)  =  &  \  \b{\lambda}(k), \qquad  k=0,\dots, N-1 \label{FOC_1}\\
 L^{\top} (\overline{\b{\nu}}_i - \underline{\b{\nu}}_i) = &   \ \b{\lambda}_i  , \hspace{0.45in} i=1, \dots, n \label{FOC_2}
\end{align}
where we have defined $\b{\lambda}_i = \bmat{\lambda_i(0),\dots, \lambda_i(N-1)}^{\top}$.
The \emph{complementary slackness condition} is given by:
\begin{align} \setlength{\itemsep}{.3in}
\b{\mu}(k) \circ \left( H ( \b{v}(k) +\b{u}(k) ) - \ccap \right) = 0,    & \qquad k=0, \dots, N-1  \label{FOC_5} \\
 \underline{\b{\nu}}_i \circ L \b{u}_{i} = 0, &    \qquad  i=1, \dots, n \label{FOC_6} \\
 \overline{\b{\nu}}_i \circ \left( \b{b}_i  - L \b{u}_i \right) = 0, &   \qquad  i=1, \dots, n.   \label{FOC_7}
\end{align}
It will occasionally prove convenient to work with alternative arrangements of the Lagrange multipliers defined above. Accordingly, we  define the vector $\b{\mu}_{\ell} = \bmat{\mu_{\ell}(0), \dots, \mu_{\ell}(N-1)}^\top$ as the sequence of Lagrange multipliers associated with each transmission line constraint $\ell = 1, \dots, 2m$. In addition, the vectors $ \underline{\b{\nu}}(k) = \bmat{ \underline{\nu}_1(k), \dots,  \underline{\nu}_n(k)}^\top $ and $\overline{\b{\nu}}(k) = \bmat{ \overline{\nu}_1(k), \dots,  \overline{\nu}_n(k)}^\top$ denote the collection of  Lagrange multipliers associated with the lower and upper bounds on storage capacity, respectively, for each time $k = 0, \dots, N-1$.

\section{Financial Storage Rights} \label{sec_FSR}

In this section, we outline the concept of \emph{financial storage rights} (FSRs), and develop their basic properties within the context of a nodal spot market for energy. Broadly, FSRs amount to financial instruments, which enable the decoupling of the ownership of storage capacity from its physical operation.
This is accomplished through the allocation of  financial property rights to storage in the form of entitlements to the merchandising surplus generated by the centralized dispatch of the storage assets.
Specifically, a FSR is defined as a sequence of hourly injections/withdrawals at a specific node in the power network, which yields  the holder a payoff according to the corresponding sequence of nodal spot prices.
Being defined as such, FSRs provide energy market participants the  ability to hedge their intertemporal exposure to  hourly price variability at specific nodes in the power network. 
And while FSRs are essentially strips of forward energy contracts, what makes this class of financial instruments unique is the fact that FSRs are issued under the physical cover of storage capacity and are funded by the surplus (\ie the intertemporal arbitrage value) that centrally operated storage generates in the spot market. 

In what follows, we investigate the \emph{revenue adequacy} of such instruments in the context of  electricity markets employing locational marginal pricing.
In particular, we  establish  conditions under which the  allocation of both financial storage and transmission rights is guaranteed to be revenue adequate, i.e., the merchandising surplus collected by the ISO is sufficient to cover the net settlement to all holders of financial storage and transmission rights.
We begin with a definition of locational marginal prices under multi-period economic dispatch with energy storage, in the following section.

\subsection{Locational Marginal Pricing}
We refer to $\b{\lambda}(k) \in \Rset^n$ as the vector of nodal prices at time $k$. More specifically, the $i^{th}$ element, $\lambda_i(k)$, denotes the price at which energy is transacted at node $i$ and time $k$. We denote by $\Lambda = \bmat{\b{\lambda}(0), \dots, \b{\lambda}(N-1)}$ the corresponding sequence of nodal prices from time $k=0$ to $N-1$. We have the following standard definitions of \emph{market equilibrium} and \emph{efficiency}. 

\begin{definition}
\label{def:effeq}
The triple $(V, U, \Lambda)$ constitutes a \emph{market equilibrium} if it satisfies \eqref{constraints_flow}, \eqref{constraints_stg} and \eqref{FOC_1}. The triple $(V, U, \Lambda)$ is said to be an \emph{efficient market equilibrium} if $(V,U)$ is also an economic dispatch.
\end{definition}

The requirement that $(V,U)$  satisfy \eqref{constraints_flow} and \eqref{constraints_stg} in Definition \ref{def:effeq} can be interpreted as \emph{market clearing} and \emph{feasibility conditions}, respectively, as they require that supply equal demand at each time period, while ensuring that the line flow and storage capacity constraints are met. Condition \eqref{FOC_1} is tantamount to requiring \emph{consumer and supplier equilibrium} at every node and time period.  In other words, relation \eqref{FOC_1} requires that the marginal cost of supply (benefit of demand) equal the nodal price $\lambda_i(k)$ for all nodes $i$ and time periods $k$. Consequently, at equilibrium, there is no opportunity for the profitable trading of energy across nodes or time.
 
It is important to note that there may exist multiple market equilibria -- some of which may not be efficient. In other words, the system operating point at a market equilibrium may not maximize social welfare. One can, however, implement an economic dispatch $(V,U)$ at a market equilibrium $(V,U, \Lambda)$, if the nodal prices $\Lambda$ are set according to \eqref{eq:LMP} -- the Lagrange multipliers derived at the  corresponding economic dispatch. Such approach to spot pricing is generally referred to as \emph{locational marginal pricing} (LMP) \cite{Schweppeetal88}.

\subsection{Merchandising Surplus} \label{subsec_MS}
 In selecting and implementing a  market equilibrium $(\V,\U,\Lambda)$,
the ISO collects payment from the consumers and remunerates the suppliers according to their respective operating points and nodal prices. In doing so, the ISO may collect a nonzero surplus. We refer to this excess as the \emph{merchandising surplus} (MS).
Indeed, it is a straightforward generalization of \cite{Wuetal96} to show that the MS can be either positive or negative at an arbitrary market equilibrium. The latter outcome is undesirable, as it may require the ISO to incur a fiscal deficit in clearing the market.
In what follows, we briefly discuss the effects of dispatch efficiency and congestion, in both transmission and storage, on the MS. First, we have a definition.

\begin{definition}
\label{def:MSdef}
The \emph{merchandising surplus} (MS) at a market equilibrium ($\V$, $\U$, $\Lambda$) is defined as
\begin{equation} \label{ms}
{\rm MS} = -\sum_{k=0}^{N-1} \sum_{i=1}^n  \lambda_i(k) v_i(k),
\end{equation}
or, equivalently, as ${\rm MS} = - \text{trace}(\Lambda^\top V)$.
\end{definition}

One can massage the expression for the MS  in \eqref{def:MSdef} to reveal the specific impact that both  transmission and storage congestion have on its value. In order to do so, we must first specify the line flows induced by the net injection profile for each time period. More formally, let $(V,U)$ be an arbitrary feasible dispatch. And denote by $p_{ij}(k)$ the resulting power flow over the line from  node $i$ to node $j$ at time $k$.\footnote{According to the formulation of DC power flow considered in Section \ref{sec:network_mod}, the line flow $p_{ij}(k)$ corresponds to a single entry of the vector $H(\b{v}(k) + \b{u}(k))$. And if there is no line connecting nodes $i$ and $j$, then $p_{ij}(k) = - p_{ji}(k) = 0$ necessarily.}
We adopt a sign convention such that  $p_{ij}(k) = - p_{ji}(k) > 0$,  if power flows from node $i$ to $j$. It follows from Kirchhoff's Current Law that \
$
v_i(k) + u_i(k)  = \sum_{j=1}^n p_{ij}(k)
$
\ for all nodes $i=1,\dots,n$. Using this relation, one can decompose the merchandising surplus as 
\begin{align} \label{eq:MS_TCS_SCS}
\ms \ = \  \tcs  \ +  \ \scs.
\end{align}
The first term in the decomposition is commonly referred to as the  \textit{transmission congestion surplus} (TCS). The second term, we refer to as the  \textit{storage congestion surplus} (SCS). Each term  satisfies:
\begin{align}
{\rm TCS} &  \ = \ \frac{1}{2} \sum_{k=0}^{N-1} \sum_{i,j=1}^n \left(\lambda_j(k) - \lambda_i(k) \right) p_{ij}(k), \label{eq:tcs} \\
{\rm SCS} &  \ = \  \sum_{k=0}^{N-1} \sum_{i=1}^n \lambda_{i}(k) u_i(k). \label{eq:scs}
\end{align}

\begin{lemma}
\label{lemma:MS}
The {\rm MS}, {\rm TCS}, and {\rm SCS} derived at an efficient market equilibrium $(\V,\U,\Lambda)$ are nonnegative quantities.
\end{lemma}

Lemma \ref{lemma:MS}  reveals an important property. Namely, at an efficient market equilibrium, the collective transactions between supply and demand are guaranteed to be \emph{revenue adequate}, i.e.,  $\ms \geq 0$.
Moreover, the reformulation of the MS in \eqref{eq:MS_TCS_SCS} reveals a decomposition of the effects due to congestion in transmission and storage on the rent collected by the ISO.

\begin{assumption}
\label{as:1}
For the remainder of the paper, we let $(V,U,\Lambda)$ denote an \emph{efficient market equilibrium}, unless otherwise specified.
\end{assumption}

\subsection{Financial Transmission Rights}

In the event that there is transmission congestion at an economic dispatch and the ISO does indeed collect a positive merchandising surplus, it is common practice in US electricity markets to reallocate the MS in the form of \emph{financial transmission rights} \cite{Hogan92, ONeillEtAl13}. Financial transmission rights can be specified in variety of ways, with the two most predominant types being defined as \emph{point-to-point} and \emph{flow-based} rights.

A \textit{point-to-point financial transmission right} (FTR) is specified in terms of  a quantity of power, a point of injection, and a point of withdrawal. It yields the holder the entitlement to receive, or obligation to pay,  the difference in nodal spot prices between the chosen point of withdrawal and the point of delivery, times the nominated quantity of power. Accordingly, an FTR may amount to a credit or liability. We have the following definition. 

\begin{definition}
A \textit{point-to-point financial transmission right} (FTR) is any triple $(i,j, \b{t}_{ij})$, where $i$ denotes an injection node, $j$ a withdrawal node, and $\b{t}_{ij} \in \Rset^N_+$ an hourly power profile spanning $N$ time periods.
The FTR  yields the holder a rent (or liability) equal to \ $(\b{\lambda}_j - \b{\lambda}_i)^{\top} \b{t}_{ij}$.
We refer to the FTR more compactly as $\b{t}_{ij}$.
\end{definition}
\begin{remark} \label{remark:FTRsymmetry}
We have implicitly required injection/extraction symmetry in our definition of FTRs, as we have considered a lossless model of power flow.
We refer the reader to \cite{PhilpottPritchard04} for a more general characterization of FTRs that accommodates lossy transmission networks.
\end{remark}

FTRs have become an important component of LMP-based electricity markets, in part, because of their ability to provide market participants with
an effective  hedge against transmission congestion costs for long-term energy transactions involving known injection and withdrawal points within the transmission network.\footnote{
We refer the reader to \cite{RosellonKristiansen13} for a recent survey on financial transmission rights.}
\textit{Flow-based} or \textit{flowgate rights} (FGRs) have also been proposed as a viable alternative or complement to FTRs \cite{ChaoPeck96, Stoft98, Chaoetal00}. 
Specifically, a FGR is a link-based transmission right, specified in terms of directed transmission link, and quantity of power flow along that link. It yields the holder the entitlement to receive a payment equal to the Lagrange multiplier (\ie shadow price) associated with the chosen link's capacity constraint multiplied the nominated quantity of power flow.  Note that the rent due to a FGR is guaranteed to be nonnegative, as the corresponding shadow prices on transmission constraints are necessarily nonnegative. 
We have the following definition of FGRs according to the model considered in this paper.

\begin{definition}
A \textit{flowgate right} (FGR) is any double $(\ell,\b{f}_\ell)$, where the  index $\ell \in \{1,\ldots,2m\}$ denotes a directed transmission link, and $\b{f}_\ell \in \Rmbb^{N}_{+}$ a hourly power profile spanning $N$ time periods.
The FGR yields the holder a rent of $\b{\mu}_\ell^{\top} \b{f}_\ell$.
We refer to the FGR more compactly as $\b{f}_{\ell}$.
\end{definition}

Although FGRs are not currently offered in the majority of transmission rights auctions that are in operation today,
the theoretical literature on the subject has converged on the viewpoint that both FTRs and FGRs could and should coexist, thereby allowing market participants the ability to decide as to what mix of rights is best   \cite{Chaoetal00, ONeillEtAl02,ONeillEtAl13}.
We adopt this perspective, and develop our mathematical results in a framework that is general enough to accommodate both types of financial rights.
Accordingly, we denote an arbitrary \emph{collection of FTRs and FGRs} by the pair $(\Tcal, \Fcal)$, where
\begin{align*}
\Tcal = \{ \b{t}_{ij} \ \left| \  i,j  = 1,\dots,n \right.\} \quad \text{and} \quad \Fcal = \{ \b{f}_\ell  \  \left| \  \ell=1,\ldots,2m \right.\}.
\end{align*} 
Here, $\b{t}_{ij}$ is the sum of all FTRs of the same type $(i,j)$, and $\b{f}_\ell$ is the  sum of all FGRs of the same type $\ell$.

In what follows, we investigate the revenue adequacy of transmission rights in nodal spot markets based on multi-period economic dispatch with storage. In particular, we establish conditions on the joint allocation of FTRs and FGRs, under which the merchandising surplus collected by the ISO is sufficient to cover their net settlement.

\begin{definition}
\label{def:FTR_rent}
The \emph{rent due to a collection of transmission rights} $(\Tcal, \Fcal)$ is defined as
\begin{align} \label{eq:FTR_rent}
\RT =  \sum_{i,j=1}^n (\b{\lambda}_j - \b{\lambda}_i)^{\top} \b{t}_{ij} \ + \ \sum_{\ell=1}^{2m} \b{\mu}_\ell^{\top} \b{f}_\ell.
\end{align}
\end{definition}
In general, the merchandising surplus collected by the ISO will not equal the rent due to a collection of transmission rights.
Their allocation must, therefore, be restricted in such a manner as to guarantee that the ISO does not incur a financial  shortfall.
A well known requirement is the \emph{simultaneous feasibility test} (SFT) \cite{Hogan92, Wuetal96, PhilpottPritchard04}.
We now extend this notion to the multi-period setting to accommodate the enlargement of the set of feasible power injections due to the  presence of storage capacity. 

We first require  additional notation.
For each time period $k = 0, \dots, N-1$,  denote by $\b{t}(k) \in \Rmbb^n$ the net injection vector induced by a collection of FTRs  in $\Tcal$, and by $\b{f}(k) \in \Rmbb^{2m}$ the vector of direction-specific flowgates induced by a collection of FGRs in $\Fcal$.\footnote{The $i$th element of the net injection vector $\b{t}(k)$ is given by $t_i(k) = \sum_{j=1}^n (t_{ij}(k)-t_{ji}(k))$. It follows that $\b{1}^{\top} \b{t}(k)=0$. Also,  the $\ell$th element of the flowgate vector $\b{f}(k)$ is given by the $k$th element of the FGR $\b{f}_\ell$. }
\begin{definition}
\label{def:sf}
A collection of transmission rights $(\Tcal, \Fcal)$ are said  to be \emph{simultaneously feasible} if there 
exists a sequence of storage injections $Q \in \Rset^{n \times N}$, which is feasible according to 
\begin{alignat*}{3}   
  \b{t}(k)  + \b{q}(k)  & \in  \  \Pmsc(\b{c} - \b{f}(k)),  \qquad &&  k=0,\dots,N-1    \\
					 	 \b{q}_{i}  &\in  \  \Umsc( \b{b}_i),  && i =1,\dots,n.   
\end{alignat*} 
\end{definition}

In other words, a collection of transmission rights $(\Tcal, \Fcal)$ are simultaneously feasible if the sequence of nodal injections induced by the FTRs in  $\Tcal$ can be reshaped by a feasible sequence of  storage injections so that the resulting nodal injections induce power flows that respect the transmission capacity limits,  derated according to the  FGRs in $\Fcal$. 
We have the following result, which establishes revenue adequacy for any simultaneously feasible collection of transmission rights.

\begin{lemma}
\label{lemma:FTR}
If $(\Tcal, \Fcal)$ are a simultaneously feasible collection of transmission rights, then their corresponding rent satisfies
\begin{align}
\RT \ \leq \ \tcs. \label{eq:FTRadequacy}
\end{align}
This inequality is \emph{tight}, in the sense that there exists a simultaneously feasible collection of transmission rights  with a corresponding rent equal to the $\tcs$.
\end{lemma}

Lemma \ref{lemma:FTR}  reveals that the transmission congestion surplus (TCS) is sufficient to cover the rent due to any collection of simultaneously feasible  transmission rights. In the event that storage facilities congests at an economic dispatch, the ISO will also collect, in addition to the TCS, a storage congestion surplus (SCS) as part of its merchandising surplus. We, therefore, define  a new class of financial instruments, which play a complementary role to transmission rights, and rely on the SCS as their primary funding source. We refer to these new instruments as  \emph{financial storage rights}.

\subsection{Financial Storage Rights} \label{sec:fsr}

\noindent We begin with a definition of \emph{financial storage rights}.

\begin{definition}
\label{def:fsr}
A \textit{financial storage right} (FSR) is any double $(i,\b{s}_i)$, where $i$ denotes a withdrawal node, and $\b{s}_{i} \in \Rset^N$ a hourly power profile spanning $N$ time periods.
The FSR  yields the holder a rent (or liability) equal to $\b{\lambda}_i^{\top} \b{s}_{i}$.
We refer to the FSR more compactly as $\b{s}_{i}$.
\end{definition}

Before embarking upon a formal analysis of FSRs and their properties, we provide a brief qualitative discussion surrounding their structure and potential use.
First, FSRs can be thought as financial property rights to storage capacity; or, more accurately, as entitlements to the intertemporal arbitrage gains that storage generates under its socially optimal operation, \ie the storage congestion surplus (SCS).
Being defined as such, FSRs enable the complete decoupling between the dispatch of storage facilities and the settlement of storage congestion charges.
Second, as tradable property rights, FSRs can be sold in forward auctions coordinated by the ISO; and the revenue generated by such auctions could serve to incentivize merchant investment in storage --   not unlike the role of FTRs in partially supporting the remuneration of merchant transmission investments \mbox{\cite{Hogan02b, KristiansenRosellon06}}.
Third, from the perspective of its holder, a FSR is equivalent to a strip of forward energy contracts.
Accordingly, FSRs yield market participants the ability to perfectly hedge physical or financial positions in the spot market against intertemporal price risk.
Such hedging capabilities represent a natural complement to FTRs and their ability to hedge spatial price risk across the network.
Finally,  an important factor  distinguishing  FSRs from standard forward energy contracts, is the crucial fact that FSRs are issued under the physical cover of storage capacity and settled against the SCS collected by the ISO, as opposed to the revenue generated from contract sales.

Different forms of financial entitlements to the storage infrastructure can be envisioned. For instance, Taylor \cite{Taylor14b} proposes an alternative form of financial storage rights, which are defined in terms of specific storage facilities, and entitle their holder to receive the shadow price on a  storage facility's  energy capacity constraint times the nominated quantity of energy. We refer to this alternative form of financial rights as \emph{energy capacity rights} (ECRs). 
Working within the confines of our idealized storage model, ECRs can be formally defined as follows.

\begin{definition}
\label{def:c-fsr}
An \textit{energy capacity right} (ECR) is any double $(i,\b{e}_i)$, where the index $i$ denotes a storage asset, and $\b{e}_i \in \Rmbb^N_+$ a hourly energy profile spanning $N$ time periods.
The ECR yields the holder a rent of  $\overline{\b{\nu}}_i^{\top} \b{e}_i$.
We refer to the ECR more compactly as $\b{e}_{i}$.
\end{definition}

In the presence of additional constraints, which limit the rate at which a storage facility can be charged or discharged, one can expand the definition of ECRs to include another class of financial rights that entitle their holder to receive the shadow price on the storage facility's power capacity constraint times the nominated quantity of power. Taylor \cite{Taylor14b} refers to such instruments as \emph{power capacity rights} (PCRs).

Definition \ref{def:c-fsr} is in contrast to our \emph{profile-based} definition of FSRs (cf. Definition \ref{def:fsr}).
Intuitively, the relationship between FSRs and ECRs is analogous to the relationship between FTRs and FGRs. And, to a large extent, the advantages and disadvantages of FSRs versus ECRs mirror those of FTRs as compared to FGRs.\footnote{We refer the reader to \cite{Chaoetal00, ONeillEtAl02,ONeillEtAl13, Oren13, Ruff01} for detailed discussions surrounding such comparisons in the context of transmission rights.} For example, while FTRs (FSRs) are convenient instruments for hedging spatial (intertemporal) price risk,  FGRs (ECRs) are instruments better suited for remunerating property rights to specific transmission lines (storage facilities).

In Section \ref{sec:generalSimultaneousFeasibility}, we present conditions on the joint offering of transmission and storage rights under which the ISO is guaranteed to be revenue adequate. To that end, we first define the rent due to a collection of FSRs and ECRs.  We denote an arbitrary \emph{collection of FSRs and ECRs} by the pair $(\Scal, \Ecal)$, where
\begin{align*}
\Scal = \{ \b{s}_i \ \left| \  i = 1,\dots,n \right.\} \quad \text{and} \quad \Ecal = \{ \b{e}_i  \  \left| \  i=1,\ldots,n \right.\}.
\end{align*} 
Here, $\b{s}_{i}$ is the sum of all FSRs of the same type $i$, and $\b{e}_i$ is the  sum of all ECRs of the same type $i$.

\begin{definition}
\label{def:rent_FSR}
The \emph{rent due to a collection of storage rights} $(\Scal, \Ecal)$ is defined as
$$ \RS = \sum_{i=1}^n \b{\lambda}_i^{\top} \b{s}_{i}  \ + \ \overline{\b{\nu}}_i^{\top} \b{e}_i. $$
\end{definition}

\subsection{A Generalized Simultaneous Feasibility Test} \label{sec:generalSimultaneousFeasibility}

We now extend our definition of multi-period simultaneous feasibility to accommodate a combination of both transmission and storage rights.
\begin{definition}
\label{def:sfGeneral}
A  collection of transmission  and storage rights $(\Tcal, \Fcal, \Scal, \Ecal)$ are said  to be \emph{simultaneously feasible} if there 
exists a sequence of storage injections $Q \in \Rset^{n \times N}$, which is feasible according to 
\begin{alignat*}{3}   
  \b{t}(k) - \b{s}(k)  + \b{q}(k)  & \in  \  \Pmsc(\b{c} - \b{f}(k)),  \qquad &&  k=0,\dots,N-1    \\
					 	 \b{q}_{i}  &\in  \  \Umsc( \b{b}_i - \b{e}_i),  && i =1,\dots,n.   
\end{alignat*} 
\end{definition}

\begin{remark}[Accommodating Inefficiencies in Storage]
\label{remark:nonidealities} While we have thus far operated under the assumption of perfectly efficient storage facilities, it is straightforward to extend the definition of simultaneous feasibility in Definition \ref{def:sfGeneral} to accommodate dissipative losses and conversion inefficiencies in storage  
by simply refining the underlying storage constraints on which it is based.
\end{remark}

Essentially, a collection of transmission and storage rights are simultaneously feasible if the nodal injections induced by the FTRs in $\Tcal$  and FSRs in $\Scal$ can be reshaped by a  sequence of storage injections, which  both respect the storage capacity constraints (derated according to the ECRs in $\Ecal$), and result in  power flows that do not violate the transmission capacity constraints (derated according to the FGRs in $\Fcal$). Notice that, in the absence of storage rights, this generalized definition of simultaneously feasibility reduces to  Definition \ref{def:sf}.
The following result characterizes the maximum rent achievable by any simultaneously feasible collection of transmission and storage rights. 

\begin{theorem}
\label{thm:MS}

If $(\Tcal, \Fcal, \Scal, \Ecal)$ are a simultaneously feasible collection of transmission and storage rights, then their corresponding rent satisfies
\begin{align}
\RT \ + \ \RS \ \leq \ \ms. \label{eq:proof_ms}
\end{align}
Moreover, this inequality is tight, in the sense that there exists a simultaneously feasible collection of rights $(\Tcal, \Fcal, \Scal, \Ecal)$  with an associated rent equal to $\ms$.
\end{theorem}

Theorem \ref{thm:MS} is reassuring, as it guarantees revenue adequacy on behalf of the ISO when jointly issuing  transmission and storage rights in a manner that is simultaneously feasible. More precisely, given a \textit{fixed} configuration of transmission and storage facilities, the MS collected by the ISO in the spot market suffices to cover  the rents of all outstanding transmission and storage rights. Revenue adequacy is not, however, guaranteed in the event of unplanned contingencies, as  the configuration of  transmission and/or storage facilities may deviate from what was assumed in the simultaneous feasibility test (SFT). 
The ISO must, therefore, specify a mechanism to compensate potential revenue shortfalls that might arise in the event that such contingencies occur.\footnote{
We refer the reader to \cite[Sec. 3.6]{Oren13}, which examines several mechanisms to cover revenue shortfalls that might occur when settling payments to FTR holders in the event of transmission line contingencies. For example, PJM   handles  revenue inadequacy in settling FTR payments  by prorating the revenue shortfall among the FTR holders; whereas, in NYISO-run markets, transmission line owners are held responsible for the shortfall \cite{ONeillEtAl13}. A mechanism of the former type generally transfers the risk of shortfall to the FTR holders, undermines the ability of FTRs to provide perfect price hedges, and is vulnerable to gaming due to the socialization of the shortfalls.
Conversely, a mechanism of the latter type fully funds the outstanding rights, thereby transferring the risk of shortfall to the transmission line owners themselves. An argument in favor of such a mechanism is that it provides an incentive to transmission line owners to effectively maintain their assets, and avoids the socialization of revenue shortfalls among the FTR holders \cite{ONeillEtAl13, Oren13}.}

\begin{remark}[Expanding the SFT via Short Positions] In \cite{Oren13}, Oren describes how  the allowance of short positions on long term flowgate rights (FGRs) might serve to incentivize the  maintenance and  incremental upgrade  of certain transmission lines to enhance their capacity in real time.  For instance, a 1 MW short FGR position on a particular line serves to enlarge the capacity of that same line by 1 MW in the simultaneous feasibility test (SFT) used to clear the forward transmission rights auction. The holder of the short FGR is paid the  forward shadow price on that line determined by the transmission rights  auction under the enlarged SFT, and is required to pay to the spot shadow price on that line in real time. Ultimately, the holder of the short FGR makes a profit if the transmission line in question turns out to be uncongested in real time. Such instruments are most naturally utilized by those market participants with the ability to maintain or upgrade transmission line ratings in real time, e.g., transmission line owners. One can envisage  an analogous role that could be played by short positions on energy capacity rights (ECRs). That is to say, the allowance of short positions on long term ECRs might serve to finance investments in storage capacity, or related technologies with load flattening capabilities, while allowing the ISO to allocate long term FSRs against such investments. Such instruments would be well suited to market participants with the ability to shape demand over time, such as load serving entities, demand response providers, and energy storage owners.
\end{remark}

\section{An Illustration of the Use of FSRs} \label{sec_example}

The natural variation of nodal spot prices over both location and time exposes market participants to price risk. 
Extreme price volatility is particularly problematic for load serving entities (LSEs) that sell electricity to  end-use customers at a fixed and predetermined price, as they face the risk that the spot price at which they pay for energy may considerably exceed the fixed price at which they are remunerated during certain hours of the day. Accordingly, LSEs and, more generally, those market participants seeking price stability in their transactions, may wish to hedge their exposure to such price risk. In the following discussion, we explain  how FSRs, in combination with 
contracts for differences (CFDs) and FTRs, can be employed to  fully hedge a long-term bilateral contract for energy, when the seller and buyer may exhibit differing intertemporal supply and demand characteristics, respectively. As a special case, the framework considered accommodates the setting in which the seller is physically constrained to deliver the contracted amount of energy through a constant power profile over a predetermined interval of time, while the  buyer is compelled to consume that amount of energy according a power profile that (predictably) fluctuates over that same interval of time, due in part to the inelastic nature of its demand.

We consider the setting in which a demander at node $j$ would like to buy an amount of energy $q_c$ (MWh) from a supplier at node $i$ to be delivered over $N$ time periods at a fixed price $\lambda_c$ (\$/MWh). The supplier is assumed to deliver this quantity of energy according to the production profile  $\b{q}_i \in \Rmbb^N$, while the demander is assumed to consume that same amount of energy according to the consumption profile $\b{q}_j \in \Rmbb^N$. While the production  and consumption profiles need not agree  at any given time period, they must balance over time, \ie $\b{1}^{\top}\b{q}_i=\b{1}^{\top}\b{q}_j=q_c$. In addition, it is assumed  that both  the demander and supplier are required to trade with the ISO  according to their respective nodal spot prices. Accordingly, the demander  pays $\b{\lambda}_j^\top\b{q}_j$, and supplier  is paid $\b{\lambda}_i^\top\b{q}_i$ in the spot market. 
Because nodal spot  prices will naturally vary over both time and location, and will therefore differ from the contract price $\lambda_c$, a hedge is required in order to execute the fixed price contract between the supplier and demander.
In what follows, we explain how a combination of a CFD, FTR, and FSR can be employed to perfectly hedge such  price differences. 

In general, the supplier will lose an amount $\lambda_c q_c - \b{\lambda}_i^\top\b{q}_i$ and the demander will gain an amount $\lambda_c q_c - \b{\lambda}_j^\top\b{q}_j$,  as a result of their respective spot market transactions.\footnote{Clearly, each of these amounts is as equally likely to be negative as positive, depending on the specific values of the contract price and nodal spot prices.} In the event that nodal spot prices are constant over both location and time, it is straightforward to see that the amount lost by the supplier is equal to the amount gained by the demander. 
Thus, a simple money transfer between the two parties in that amount is sufficient to perfectly hedge the fixed price contract. Such transfer can be accomplished with a CFD, which requires that the demander pay the supplier the amount $\lambda_c q_c - \b{\lambda}_i^\top\b{q}_i$, as is recorded in the second line of Table \ref{tab:example}.

More generally, the CFD specified above leaves the demander exposed  to a \emph{transmission congestion charge} when nodal spot prices vary over location,  and  a \emph{storage congestion charge} when nodal spot prices vary over time.\footnote{Of course, this is but one of several natural ways in which the CFD might be specified. Alternative specifications that entail risk sharing between the supplier and demander can also be envisaged.}
The specific form these congestion charges is made explicit in the following expression \eqref{eq:exposure}, which disentangles the individual effects that locational and temporal price differences have on the demander's market exposure after having settled the CFD with the supplier.
\begin{align} \label{eq:exposure}
\underbrace{\b{\lambda}_j^\top\b{q}_j}_{\substack{\text{spot market} \\ \text{charge}}} \hspace{-.05in}  + \  \underbrace{(\lambda_c q_c - \b{\lambda}_i^{\top} \b{q}_i)}_{\text{CFD charge}}  \ =  \   \lambda_c q_c   \  + \hspace{-.05in}  \underbrace{(\b{\lambda}_j-\b{\lambda}_i)^{\top} \b{q}_i}_{\substack{\text{transmission congestion} \\ \text{charge}}} + \  \underbrace{\b{\lambda}_j^\top (\b{q}_j - \b{q}_i)}_{\substack{\text{storage congestion} \\ \text{charge}}}  
\end{align}
It is straightforward to see that the transmission (storage) congestion charge vanishes when the nodal spot prices are constant across location (time). 
 In the event that nodal spot prices vary over location, the resulting  transmission congestion charge can be perfectly hedged with a FTR from node $i$ to node $j$ given by $\b{t}_{ij} := \b{q}_i$. Similarly,  a FSR at node $j$ given by $\b{s}_j := \b{q}_j - \b{q}_i$  yields a perfect hedge against the storage congestion charge, in the event that nodal spot prices vary across time.\footnote{It is worth mentioning that, should the two parties enter into a bilateral contract specifying common production and consumption profiles, \ie $\b{q}_i = \b{q}_j$, the storage congestion charge would vanish -- thereby eliminating the  need for the procurement of a FSR in the pursuit of a perfectly hedged, fixed price contract.} Essentially, this FSR yields the demander a hedge, which is identical to that which could have been produced using a physical storage facility  to purchase the profile $\b{q}_i - \b{q}_j$ in the spot market at node $j$.

In combination with the CFD,  the procurement of a FTR and a FSR by the demander yields a perfectly hedged, fixed price contract between the supplier and demander -- provided that both  parties deliver and consume power according to the profiles specified by the contract. However, as is argued in  \cite{BushnellStoft97}, such price risk cannot be costlessly eliminated, as FTRs and FSRs will, in general, have nonzero value in expectation. 
We refer the reader to Table  \ref{tab:example} for a detailed accounting of the transactions described in the preceding discussion.

\setlength\extrarowheight{4pt}

\begin{table}
\caption{Using CFDs, FTRs, and FSRs to synthesize bilateral contracts.}
\begin{center}
\begin{tabular*}{\textwidth}{>{\centering\arraybackslash}m{0pt}@{}llllll}
\hline  
	&& \multirow{2}{*}{Contract or Market}		&	\multicolumn{2}{l}{Supplier at node $i$} &	\multicolumn{2}{l}{Demander at node $j$}	\\\cline{4-7}
	&& 				&	Quantity	&	Payment							&	Quantity	&	Payment		\\\hline
&1	& Spot market	&	$\b{q}_i$	&	$\b{\lambda}_i^{\top} \b{q}_i$	&	$-\b{q}_j$	&	$-\b{\lambda}_j^{\top} \b{q}_j$ 	\\
&2	& CFD 			&	$\b{q}_i$	&	$\lambda_c q_c - \b{\lambda}_i^{\top} \b{q}_i$	&	$-\b{q}_i$	&	$-(\lambda_c q_c - \b{\lambda}_i^{\top} \b{q}_i)$	\\
&3	& FTR 						&	--			&	--	&	$\b{t}_{ij}:=\b{q}_i$		&	$(\b{\lambda}_j-\b{\lambda}_i)^{\top} \b{t}_{ij}$\\
&4	& FSR 						&	--			&	--	&	$\b{s}_j:=\b{q}_j-\b{q}_i$	&	$\b{\lambda}_j^{\top} \b{s}_j$		\\\hline
&5	& Total				&	--	&	$\lambda_c q_c$				&	--	&	$-\lambda_c q_c$	\\\hline
\end{tabular*}
\end{center}
\label{tab:example}
\end{table}

\section{Conclusion and Future Work} \label{sec_Conclusions}

In this paper, we have proposed a general regulatory and market framework to enable the \emph{open access} integration  of storage, in which storage is treated as a communal asset accessible by all market participants.
Such an approach represents a substantial departure from the more standard storage integration paradigm in which a storage owner-operator pursues her individual profit maximizing interests within the confines of her local spot market.
Central to our proposal is the concept of \emph{financial storage rights} (FSRs), which are defined as a sequence of nodal power injections and withdrawals that yield the holder a  payment according to the corresponding sequence of nodal prices.
Qualitatively, FSRs represent financial property rights to the capacity of centrally operated storage facilities.
This is in sharp contrast to the physical rights proposed in \cite{Heetal11}.
An essential advantage of FSRs and the modus operandi they entail is that their allocation does not interfere with the socially optimal operation of storage or the independence of the ISO, regardless of the ownership structure of the storage facilities.
Most importantly, FSRs enable the synthesis of fully hedged, fixed-price bilateral contracts for energy, when the seller and buyer exhibit differing intertemporal supply and demand characteristics, respectively.

More broadly, we envision storage owners trading such FSRs with other market participants through (short and long term) forward auctions and secondary markets centrally coordinated by the ISO; not unlike markets for FTRs today.
And, by selling financial rights to their energy storage capacity in various forward auctions (\eg yearly, quarterly, etc.), storage owners can more finely manage their exposure to spot price volatility. In addition, the allowance of short positions on long term ECRs might serve to finance investments in storage capacity, or related technologies with load flattening capabilities, while allowing the ISO to issue long term FSRs against such investments. Also, the auction revenue derived from the forward sale of FSRs may serve as a transparent long term market signal to partially guide  merchant  investment in storage.

The study of financial storage rights presented in this paper represents  an initial point of analysis. Many interesting questions remain.
First, how should the ISO structure an auction mechanism to jointly allocate both financial transmission and storage rights?
For instance, the simultaneous feasibility test (SFT) that we propose would require coordination in clearing both the transmission and storage right auctions.
This might be too cumbersome to be practical. Accordingly, it would be of interest to explore the design of alternative conditions for simultaneous feasibility that would enable the decoupling of transmission and storage right auctions.
Second, the potential value that energy storage offers to the power system goes well beyond the application of intertemporal energy arbitrage considered in this paper \cite{Sioshansietal12}.
For example, certain storage technologies posses the capability of providing voltage support or frequency regulation services.
A natural question then, is how might one expand the concept of FSRs to incorporate these value streams as well?  
Third, it would be of interest to generalize the market framework considered to accommodate a broader family of technologies capable of shifting energy in time (\eg flexible demand-side resources).

\section{Acknowledgments}
The authors would like to thank Shmuel Oren, Pravin Varaiya (UC Berkeley), Michael Swider (New York ISO),  Eugene Litvinov, Feng Zhao, Chris Geissler, Tongxin Zheng,  Jinye Zhao (New England ISO), participants of the 2014 CERTS Reliability and Markets Internal Program Review, and an anonymous referee for their helpful comments and feedback.
This work was supported in part by NSF grant ECCS-1351621, NSF grant CNS-1239178, NSF grant IIP-1632124, NSF grant CNS-1135844, NSF grant US DoE under the CERTS initiative, and the Simons Institute for the Theory of Computing.



\bibliographystyle{plain}

\appendix
\label{sec:appendix}

\section*{Proof of Lemma \ref{lemma:MS} }

Let $(\V,\U,\Lambda)$ denote an efficient market equilibrium throughout.
And, let ($\b{\lambda}(k)$, $\gamma(k)$, $\b{\mu}(k)$, $\overline{\b{\nu}}_i$,  $\underline{\b{\nu}}_i$) denote the corresponding Lagrange multipliers satisfying the KKT conditions \eqref{constraints_flow}-\eqref{FOC_7} for all $k$ and $i$. 
We prove the desired result by establishing nonnegativity of both the TCS and SCS.

\begin{proposition}
\label{prop:TCS}
$ \tcs = \sum_{k=0}^{N-1} \b{\mu}(k)^{\top} \b{c} \ge 0$.
\end{proposition}

\begin{proof}
We have that $\tcs = -\sum_{k=0}^{N-1} \b{\lambda}(k)^{\top} (\b{v}(k)+\b{u}(k))$ based on its definition in \eqref{eq:tcs}.
Substituting for $\b{\lambda}(k)$ according to Equation \eqref{eq:LMP}, and using the fact that $\b{1}^{\top}(\b{v}(k)+\b{u}(k)) = 0$ for all $k$, we have that  
\begin{align*}
\tcs = \sum_{k=0}^{N-1} \b{\mu}(k)^{\top} H (\b{v}(k)+\b{u}(k)).
\end{align*}
The complementary slackness condition  \eqref{FOC_5} yields $\tcs = \sum_{k=0}^{N-1} \b{\mu}(k)^{\top} \b{c}$, which is clearly nonnegative.  \hfill $\blacksquare$
\end{proof}

\begin{proposition}
\label{prop:SCS}
$\scs = \sum_{i=1}^{n} \overline{\b{\nu}}_i^{\top} \b{b}_i \ge 0$.
\end{proposition}

\begin{proof}
We have that  $\scs  = \sum_{i=1}^{n} \b{\lambda}_i^{\top} \b{u}_i$ based on its definition in \eqref{eq:scs}. A direct substitution of the stationarity condition \eqref{FOC_2} and complementary slackness conditions \eqref{FOC_6}-\eqref{FOC_7} yields 
\begin{align*}
\scs = \sum_{i=1}^{n} \overline{\b{\nu}}_i^{\top} L \b{u}_i - \underline{\b{\nu}}_i^{\top} L \b{u}_i = \sum_{i=1}^{n} \overline{\b{\nu}}_i^{\top} \b{b}_i,
\end{align*}
which is clearly nonnegative. \hfill $\blacksquare$
\end{proof}
The desired result follows from Propositions \ref{prop:TCS} and \ref{prop:SCS}, as $\ms = \tcs + \scs$.

\section*{Proof of Lemma \ref{lemma:FTR}}

Let $(V,U,\Lambda)$ denote an efficient market equilibrium throughout the proof. Also, let $\gamma(k)$ and  $\b{\mu}(k)$  denote the corresponding Lagrange multipliers  associated with the power balance and line flow capacity constraints for each time period $k=0, \dots, N-1$. 
The maximum rent achievable by any simultaneously feasible collection of transmission rights is given by the optimal value of 
\begin{align}
& \text{maximize}   &&	\RT && \label{eq:lem2:objective}\\
& \text{subject to} 			&&	\b{t}(k) + \b{q}(k)  \in  \  \Pmsc(\b{c}-\b{f}(k)), &&	k=0,\dots,N-1  \label{lem2:constraints_flow}\\
						&&& \b{q}_{i}  \in  \  \Umsc(\b{b}_i),  && i =1,\dots,n, 
\label{lem2:constraints_stq}
\end{align}
where $\RT = \sum_{k=0}^{N-1} -\b{\lambda}(k)^{\top} \b{t}(k) \ + \  \b{\mu}(k)^{\top} \b{f}(k)$.
It is not difficult to show (by induction) tha any feasible solution of problem \eqref{eq:lem2:objective}-\eqref{lem2:constraints_stq} must satisfy $\b{q}_i=\b{0}$ for  $i = 1, \dots, n$.
This stems from the injection/extraction symmetry required by our definition of FTRs (which implies that $\b{1}^{\top} \b{t}(k) = 0$ for $k= 0, \dots, N-1$), and our assumption of zero initial stored energy.
One can, therefore, equivalently reformulate problem \eqref{eq:lem2:objective}-\eqref{lem2:constraints_stq} as
\begin{align}
& \text{maximize}   &&	\RT && \label{eq:lem2:objective2}\\
& \text{subject to} 	&&	\b{t}(k)  \in  \  \Pmsc(\b{c}-\b{f}(k)), \qquad 	k=0,\dots,N-1.  & &\label{eq:lem2:constraints_flow2}
\end{align}
This is a convex optimization problem with linear constraints in the decision variables  $\b{t}(k)\in \Rmbb^n$ and $\b{f}(k) \in \Rmbb^{2m}_+$ ($k=0, \dots, N-1$).
It follows that a primal optimal solution is characterized by the existence of Lagrange multipliers $\tilde{\gamma}(k)\in\Rmbb$ and $\tilde{\b{\mu}}(k) \in \Rmbb^{2m}_+$ ($k=0, \dots, N-1)$ such that the KKT conditions \eqref{eq:lem2:constraints_flow2}-\eqref{eq:proof:FOC_3} hold.
The stationarity condition is given by:
\begin{align}
\tilde{\gamma}(k) \b{1}  - H^{\top} \tilde{\b{\mu}}(k)  - \b{\lambda}(k)= & \ 0, \qquad  k=0,\dots, N-1 \label{eq:proof:FOC_1}\\
\tilde{\b{\mu}}(k) - \b{\mu}(k) = & \ 0, \qquad  k=0,\dots, N-1. \label{eq:proof:FOC_2}
\end{align}
The complementary slackness condition is given by:
\begin{align} \setlength{\itemsep}{.3in}
\tilde{\b{\mu}}(k) \circ \left( H \b{t}(k) - \ccap + \b{f}(k) \right) = 0,    & \qquad k=0, \dots, N-1.  \label{eq:proof:FOC_3}
\end{align}
Recall that $(V,U,\Lambda)$ and  $\gamma(k)$, $\b{\mu}(k)$ satisfy the  KKT conditions \eqref{constraints_flow}, \eqref{eq:LMP} and \eqref{FOC_5} associated with multi-period economic dispatch problem. It follows that  a primal optimal solution to problem \eqref{eq:lem2:objective2}-\eqref{eq:lem2:constraints_flow2} is given by 
\begin{align*}
\b{t}(k) = \b{v}(k) + \b{u}(k) \quad \text{and} \quad \b{f}(k)=\b{0}, \qquad k=0, \dots, N-1.
\end{align*}
This optimal solution yields a collection of transmission rights with an associated  rent of
$$
\RT = \sum_{k=0}^{N-1} -\b{\lambda}(k)^{\top} \b{t}(k) = \sum_{k=0}^{N-1} -\b{\lambda}(k)^{\top} (\b{v}(k)+\b{u}(k)).
$$
Upon examination of Equations \eqref{ms}-\eqref{eq:scs}, it is straightforward to verify that the right-hand side equals the $\tcs$ associated with $(V,U,\Lambda)$, thus completing the proof.

\section*{Proof of Theorem \ref{thm:MS}}

Let $(V,U,\Lambda)$ denote an efficient market equilibrium, and  let $\gamma(k)$, $\b{\mu}(k)$ ($k=0, \dots, N-1$) and $\b{\overline{\nu}}_i$, $\b{\underline{\nu}}_i$ ($i =1, \dots, n$) denote the corresponding Lagrange multipliers.
The maximum rent achievable by any simultaneously feasible collection of transmission and storage rights is given by the optimal value of
\begin{align}
& \text{maximize}   &&	\RT + \RS && \label{eq:thm1:objective}\\
& \text{subject to} 			&&	\b{t}(k) - \b{s}(k)  + \b{q}(k)  \in  \  \Pmsc(\b{c}-\b{f}(k)), &&	k=0,\dots,N-1  \label{eq:thm1:constraints_flow}\\
						&&& \b{q}_{i}  \in  \  \Umsc(\b{b}_i-\b{e}_i), && i =1,\dots,n,
\label{eq:thm1:constraints_stq}
\end{align}
where  the objective function is given by
\begin{align*}
\RT + \RS 
&= \sum_{k=0}^{N-1} \left( -\b{\lambda}(k)^{\top} \b{t}(k) + \b{\mu}(k)^{\top} \b{f}(k) \right) \  + \ \sum_{i=1}^n \left( \b{\lambda}_i^{\top} \b{s}_i + \overline{\b{\nu}}_i^{\top} \b{e}_i \right).
\end{align*}
This is a convex optimization problem with linear constraints in the decision variables  $\b{t}(k)$, $\b{f}(k)$, $\b{q}(k)$, $\b{s}(k)$  ($k=0, \dots, N-1$)  and $\b{e}_i$ ($i =1, \dots, n$).
Using arguments identical to those employed in the proof of Lemma \ref{lemma:FTR}, it is straightforward to  verify that an optimal solution to problem \eqref{eq:thm1:objective}-\eqref{eq:thm1:constraints_stq} is given by 
\begin{align*}
\b{t}(k)=\b{v}(k)+\b{u}(k), \quad \b{f}(k)=\b{0}, \quad \b{q}(k)=\b{u}(k), \quad \b{s}(k)=\b{u}(k), \quad \text{and} \quad \b{e}_{i}=\b{0}
\end{align*}
for all $k= 0, \dots, N-1$ and $i  = 1, \dots, n$.
This optimal solution yields a collection of transmission and storage rights with an associated rent of
\begin{align*}
\RT + \RS = -\sum_{k=0}^{N-1} \b{\lambda}(k)^{\top} \b{v}(k).
\end{align*}
According to \eqref{ms}, this equals the $\ms$ associated with $(V,U,\Lambda)$,  thus completing the proof.

\end{document}